\documentstyle[12pt]{article}
\title{A family of \'etale coverings of the affine line}
\author{Kirti Joshi}
\date{}
\begin{document}
\maketitle
\newtheorem{thm}{Theorem}
\newtheorem{lemma}{Lemma}
\newtheorem{conj}{Conjecture}
\newenvironment{proof}{
                       \trivlist \item[\hskip \labelsep{\bf Proof}:]
                      }{
			\hfill$\Box$\endtrivlist
		      }
\def\C{{\bf C}}
\def\P{{\bf P}}
\let\into=\hookrightarrow
\def\A{{\bf A}}
\def\End{\mathop{\rm End}\nolimits}
\def\Ga{{\cal G}_a}
\def\alg{{\rm alg}}
\def\F{{\bf F}}
\def\Gal{\mathop{\rm Gal}\nolimits}
\def\Z{{\bf Z}}
\let\surjects=\longrightarrow
\let\isom=\simeq

\section{Introduction}
	This note was inspired by a colloquium talk given by
S.~S.~Abhyankar at the Tata Institute\footnote{During his visit in the
month of December, 1992}, on  the work of Abhyankar, Popp and Seiler
(see \cite{Popp}). It was pointed out in this talk that classical
modular curves can be used to construct (by specialization) coverings
of the affine line in positive characteristic.  In this ``modular''
optic it seemed natural to consider Drinfel'd modular curves  for
constructing coverings of the affine line.  This note is a direct
outgrowth of this idea.

	 While it is trivial to see that affine line in characteristic
zero has no non-trivial \'etale coverings, in \cite{Abhyankar} it was shown
that the situation in positive characteristic is radically different
and  far more interesting. Let us, for the sake of definiteness, work
over a field $k$ of characteristic $p>0$.  In \cite{Abhyankar},
Abhyankar  conjectured that any finite group whose order is divisble
by $p$ and which is generated by its $p$-Sylow subgroups (such a
finite group is sometimes called a ``quasi-$p$-group''), occurs as
quotient of the algebraic fundamental group of the affine line. It is
customary to write $\pi_1^\alg(\A^1_k)$ to denote the algebraic
fundamental group of the affine line.  While Abhyankar's conjecture
indicates that the algebraic fundamental group of the affine line is
quite complicated, our result  perhaps illustrates its cyclopean
proportions. The result we prove (see Theorem~\ref{main theorem}
below) is the analogue of the following well-known classical result,
which falls out of the theory of elliptic modular curves over the
field of complex numbers: there is a continous quotient
$\pi_1^\alg(\P^1_{\C}-\{0,1728,\infty\}) \surjects \prod_{p}
SL_2(\Z_p)/\{\pm 1\}$, where $\Z_p$ denotes the ring of $p$-adic
integers.

	We would like to thank N. Mohan Kumar for numerous suggestions
and conversations; he also explained to us a variant of Abhyankar's
Lemma, which is crucial to our argument. We would also like to thank
E.-U.~Gekeler, Dipendra Prasad for electronic correspondence while
this note was being written -- their remarks and suggestions have been
extremely useful; in particular we would like to point out that the
conjecture stated at the end of this note was formulated with the help
of Dipendra Prasad. We would also like to thank M.~V.~Nori for useful
comments and Dinesh Thakur for encouragement.

\section{Resum\'e of Drinfel'd modules and their moduli}
	In this section we recall a few of the standard facts about
Drinfel'd modules. Since the basic theory of Drinfel'd modules and
their moduli is well documented we will be brief; all the facts
which we will need can be found in the following standard references:
\cite{Drinfel'd}, \cite{Deligne-Husemoller}, \cite{Gekeler2}.
Since we do not need the full strength of Drinfel'd's work, we will
work with a very special situation which is required for our purpose.
In this section, we will outline this special situation.

	Let us  fix some notations. Let $\F_q$ denote a finite field
with $q=p^m$ elements and of characteristic $p$. We will write $A =
\F_q[ T ]$, $K = \F_q(T)$. Further denote by $K_\infty$, the
completion of $K$ along the valuation corresponding to $1/T$. Denote
by $C$ the completion of the algebraic closure of $K_\infty$. The
field $C$ is thus a ``universal domain'' of characteristic $p$. There
is a natural inclusion $K \into C$.

	Let $\Ga$ denote the additive group scheme.  Then it is easy
to see that the ring of endomorphisms of $\Ga$ defined over $C$,
denoted by $\End_C(\Ga)$, is a noncommutative ring generated by the
Frobenius endomorphism. More precisely, for an indeterminate $\tau$,
consider the noncommutative ring $C\{\tau_p\}$ of all polynomials in
$\tau$ with coefficients in $C$, and where the multiplication rule is
given by $\tau_p a = a^p \tau_p$ for all $a\in C$.  Then one checks that
$\End_C(\Ga) \isom C\{\tau_p\}$. Observe that this isomorphism gives
rise to a ring homomorphism $\partial:\End_C(\Ga) \to C$ which sends
an endomorphism of $\Ga$ to the constant term of the polynomial in
$\tau$  associated to it.

	A {\em Drinfel'd module} $\phi$ over $C$ is a ring
homomorphism $\phi:A \to \End_C(\Ga)$ such that composite map
$\partial \phi $ is the natural inclusion of $A\to C$. It is easy to
check that any such map factors through the subring
$C\{\tau_p^m\}\subset C\{\tau_p\}$. For simplicity of notation, we
will write $\tau=\tau_p^m$.  Thus for any $a\in A$ we have an
endomorphism $\phi_a$ of $\Ga$ defined over $C$. Moreover, note that,
any such $\phi$ is $\F_q$ linear.  So as $A$ is generated as an
$\F_q$-algebra by $T$, giving a $\phi$ thus amounts to specifying a
single endomorphism of $\Ga$ corresponding to $T\in A$,  $\phi_T =
\sum_{i=0}^r a_i \tau^i$, with $a_0 =T$. Then this can be extended to
all of $A$. The $\tau$-degree of $\phi_T$ is called the {\em rank} of
the Drinfel'd module $\phi$.

	If $\phi,\phi'$ are two Drinfel'd modules then a {\em morphism
of Drinfel'd modules} is an endomorphism $u\in\End_C(\Ga)$ such that
for all $a\in A$ we have $\phi_a u = u \phi_a' $. It is easy to see
that any such $u$ is in fact contained in $C\{\tau\}$.

	For any $a\in A$ a Drinfeld module $\phi$ specifies a closed
subgroup-scheme of $\Ga$: the kernel of the endomorphism $\phi_a$,
$\ker(\phi_a)$. For example, if $a=T$, then the kernel of $\phi_T$ is
simply the roots of the polynomial $T+\sum_{i=1}^r a_i X^{q^i-1} =0$
together with $0$.  Moreover, for any ideal $I\subset A$, we can
define a subgroup scheme of $\Ga$ using the ring structure. One checks
that $\ker(\phi_I)(C)$ is a free  $(A/I)$-module of rank $r$, where $r$
is the rank of $\phi$.  This lets us define a notion of an $I$-level
structure on $\phi$.  Drinfel'd has shown (see \cite{Drinfel'd})
that there is a moduli of Drinfel'd modules with level structure (in
general we have only a ``coarse moduli scheme''). Our main interest is
the case of rank two Drinfel'd modules. And henceforth, we shall work
with Drinfel'd modules of rank two.

	The theory of Drinfel'd modules of rank two behaves like the
theory of elliptic curves. The fact that the moduli of Drinfel'd
modules of rank two over $C$, is a smooth affine curve over $C$ is a
very special case of a fundamental result of Drinfel'd (see
\cite{Drinfel'd}). We need several facts about  these Drinfel'd
modular curves. The first fact we need is the following:

\begin{thm}
	 The ``coarse moduli'' of rank two Drinfel'd modules over $C$
is the affine line over $C$.
 \end{thm}
 \begin{proof}
	For a proof see \cite{Goss} paragraph 1.32, or
\cite{Gekeler1}.
 \end{proof}

	Thus the above situation is analogous to the classical
situation for elliptic curves. In fact the identification with the
line is given by a ``j-invariant''. If $\phi_T = T+ a\tau + b \tau^2$,
is a rank 2 module then $b\neq 0$ and then its $j$-invariant is
$a^{q+1}/b$.

Before we need some notations. For any $\F_q$ algebra $R$, let
	$$G_1(R) = \left\{g\in GL_2(R)\big| \det(g)\in \F_q^* \right\},$$
	$$ G(R) = G_1(R)/Z(\F_q),$$
 where $Z(\F_q)$ is the group of $\F_q$-valued scalar matrices with
nonzero determinant.

	The following result  is really the heart of our construction.

\begin{thm}\label{Drinfel'd-Gekeler}
	 For every non-zero ideal $I\subset A$, there is a ``coarse
moduli'' of Drinfel'd modules of rank two with full $I$-level structure
exists and is an affine curve over $C$. There is a ``forget the level
structure'' morphism to $\A^1_C$. This map is branched over $0\in
\A^1_C$. The covering is Galois and the Galois group is $G(A/I)$. The
ramification index of any point lying over $0$ is $q+1$ and is
independent of $I$.  In particular, the ramification is tame.
 \end{thm}

\begin{proof}
	As mentioned earlier, the existence of the moduli is due to
Drinfel'd (see \cite{Drinfel'd}). These curves have been studied in
great detail by Goss and Gekeler.  In the  our case the Galois group
can easily be calculated, a convenient reference for it is
\cite{Gekeler2}, for instance see Lemma 1.4 on page 79, also see the
first section of \cite{Gekeler1}.  The ramification information is
computed in \cite{Goss}, Lemma 4.2.  One also finds it computed in
\cite{Gekeler2}, on page 87.

	As in the classical situation, the ramification takes place
over the Drinfel'd module with extra automorphisms. From the
definitions, it is clear that an automorphism of a Drinfel'd module of
rank two over $C$ is firstly an automorphism $u$ of $\Ga$, defined
over $C$.  Clearly any such automorphism must be an invertible element
of $C$.  Then a simple calculation shows that if a nonzero element of
$C$ is an automorphism of a Drinfel'd module then it must be a root of
unity.  One checks that with the exception of the module with
$j$-invariant equal to zero, the automorphism group of the Drinfel'd
module is $\F_q^*$. The module with $j$-invariant equal to zero has
automorphism group $F_{q^2}^*$.  Thus in particular, the ramification
is tame. These facts are easily proved by explicit calculations.
 \end{proof}

\section{The main theorem}
	We are now ready to state and prove our main theorem. One
should note that most  of the work has already gone in the
construction and analysis the  of the moduli of Drinfel'd modules.

\begin{thm}\label{main theorem}
	 If $p=2,3$ then assume that $q=p^m, m\geq 2$.  There is a
family of \'etale coverings, $Y_I$ of $\A^1_C$, indexed by the nonzero
proper ideals $I\subset A$. The curves $Y_I$ are affine curves over
$C$ and the covering $Y_I \to \A^1_C$ is Galois with Galois group
$SL_2(A/I)/\{\pm1\}$. Moreover, these coverings form an inverse system
indexed by $I$. Thus in particular we have a continous quotient
	$$\pi_1^\alg( \A^1_C) \to \lim_{\longleftarrow \atop I}
		\left( SL_2(
			A/I )/\{\pm 1\} \right).$$
  \end{thm}

\begin{proof}
	By Theorem~\ref{Drinfel'd-Gekeler}, we have a tamely ramified
covering of the affine line which is branched over one point. Also
note that the ramification index of any point over $0\in \A^1_C$ is
$q+1$, independent of the ideal $I\subset A$.

	Now we  apply a suitable variant of Abhyankar's Lemma to
remove the tame ramification.  The crucial thing is to ensure, if
possible, that the ``pull back'' coverings remain irreducible. The
following variant of Abhyankar's Lemma (see Lemma~\ref{Mohan's
lemma}) which was pointed out to me by N. Mohan Kumar, gives an
explicit criterion to check irreducibility, then we apply this
criterion to the case at hand.  This is an easy exercise in elementary
group theory. We have stated all the necessary results as a sequence
lemmas, and since the proofs of all the individual statements are
easy, we leave the details to the reader.
 \end{proof}

\begin{lemma}\label{Mohan's lemma}
	 Let $k$ be an algebraically closed field of characteristic
$p$. Let $X \to \A^1_k$ be a finite Galois cover defined over $k$,
with Galois group $G$. Further assume that the cover is branched over
$0\in \A^1_k$, and any point lying over it is tamely ramified with
ramification index $n$. Let $\A^1_k \to \A^1_k$ be a $\Z/n$ covering
ramified completely at $0$, and unramified elsewhere. Let $X'$ be the
normalization of the fibre product $X \times_{\A^1_k} \A^1_k$. Suppose
that there are no nontrivial homomorphisms $G \to \Z/n$. Then $X'$ is
irreducible, and the Galois group of the covering $X' \to \A^1$ is
$G$.
 \end{lemma}

\begin{proof}
	Clearly one is reduced to proving the following field theory
statement:

	Let $L/k(t)$ be a finite Galois extension with Galois group
$G$. Let $E=k(t^{1/n})$. Then $L/k(t)$ and $E/k(t)$ are linearly
disjoint over $k(t)$ if and only if there are no nontrivial
homomorphisms $G \to \Z/n$.

	And the above statement is immediate from the fact that
$E/k(t)$ is a Kummer extension. This finishes the proof.
 \end{proof}

	We now need some elementary group theoretic lemmas to apply
the above criterion. Since the proofs are easy we will state the
lemmas without proofs.

\begin{lemma}\label{lemma1}
	 Let $\wp \neq 0$ is a prime ideal of $A$. Then for all $k
\geq 1$ and for all $n\geq 2$, the natural morphism
	$$GL_k(A/\wp^n) \to GL_k( A/\wp^{n-1})$$
 is surjective.
 \end{lemma}

Recall that for any $\F_q$ algebra $R$, we had defined the groups
	$$G_1(R) = \left\{g\in GL_2(R)\big| \det(g)\in \F_q^* \right\},$$
	$$ G(R) = G_1(R)/Z(\F_q),$$
 where $Z(\F_q)$ is the group of $\F_q$-valued scalar matrices with
nonzero determinant. Further, for any prime ideal $\wp\neq 0$ of $A$,
write $G_1^n = G_1(A/\wp^n)$, and $G^n = G(A/\wp^n)$, for all $n\geq
1$. The results which follow are valid for any non-zero prime ideal
$\wp$, so we have supressed $\wp$ in our notations.

\begin{lemma}\label{lemma2}
	 The natural map $G_1^n \to G_1^{n-1}$ is surjective for all
$n\geq 2$.
 \end{lemma}

\begin{lemma}\label{lemma3}
	 The natural map $G^n \to G^{n-1}$ is surjective for all
$n\geq 2$.
 \end{lemma}

	Denote by $G^{(n,n-1)}$ the kernel of the map $G^n \to
G^{n-1}$, similarly write $G_1^{(n,n-1)}$, for the kernel of the
corresponding map for $G_1$.

\begin{lemma}\label{lemma4}
	 Let $q=2^m, m\geq 2$. Let $F/\F_q$ be any finite
extension. Then there are no nontrivial morphisms $G(F) \to
\Z/(q+1)$.
 \end{lemma}

\begin{lemma} \label{lemma5}
	 Let $q=2^m,m\geq 2$. Then for any $n\geq 1$ there are no
nontrivial maps $G^n \to \Z/(q+1)$.
 \end{lemma}
 \begin{proof}
    The proof is by induction on $n$. For $n=1$, we are done by the
previous lemma. Now show that the kernel of any map $G^n \to \Z/(q+1)$
contains $G^{(n,n-1)}$. Then we are done by induction.
 \end{proof}

\begin{lemma}\label{lemma6}
	 Let $q = p^m, p\neq 2$. If $p=3$ then $m\geq 2$.
Then there is a canonical morphism $G^n \to \F_q^*/\F_q^{*2}$. In
particular we have a natural map $G^n \to \Z/2$, obtained by
identifying $\F_q^*/\F_q^{*2}$ with $\Z/2$.
 \end{lemma}

 \begin{lemma}\label{lemma7}
	 Let $q=p^m, p\neq 2$. If $p=3$ then $m\geq2$. Let $F/\F_q$ be
any finite extension. Then any morphism $G(F)\to \Z/(q+1)$ factors
through the canonical morphism given by the above lemma, followed by
the inclusion of $\Z/2\to \Z/(q+1)$.
 \end{lemma}

 \begin{lemma}\label{lemma8}
	 Let $q = p^m, p\neq 2$. If $p=3, m\geq 2$. Then any
nontrivial morphism $G^n \to \Z/(q+1)$ factors through the canonical
morphism.
 \end{lemma}

 \begin{proof}
	This is again proved by induction on $n$. As before, one
checks that any such map is trivial on $G^{(n,n-1)}$.
 \end{proof}

Now we can identify our Galois groups.

\begin{lemma}\label{lemma9}
	 Let $q=2^m, m \geq 2$. For any non zero prime ideal
$\wp\subset A$, we have:
	$$G^n = SL_2(A/\wp^n).$$
 \end{lemma}

\begin{lemma}\label{lemma10}
	 Let $q=p^m, p\neq 2$, if $p=3$ then $m\geq 2$. let
$\wp$ be any non zero prime ideal in $A$. For any $n\geq 1$, let
	$$\tilde{G^n} = ker(G^n \to \Z/2).$$
 Then $\tilde{G^n} = SL_2(A/\wp^n)/\{\pm 1\}$.
 \end{lemma}

	Thus we can now prove our main theorem. If $p=2$, then the
pull back coverings remain irreducible and the Galois group is $G^n$.
If $p$ is odd, then there is a quadratic subfield in common. And the
Galois group is  $\tilde{G^n}$. Then we are done by the above lemmas.

	After the result of Madhav Nori (see \cite{Madhav}) and the
one proved above, we would like to advance the following conjecture:
 \begin{conj}
	 Let $G/K$ be any isotropic, semisimple, simply-connected
algebraic group over $K=\F_q(T)$, with center $Z$. Let $\A^{{\rm
fin}}_K$ denote the finite adeles of $K$.  Then any maximal compact
subgroup of $G(\A^{{\rm fin}}_K)/Z(\F_q)$ occurs as a continous
quotient of the fundamental group of the affine line over $C$.
 \end{conj}


\flushleft{School of Mathematics,\\
			Tata Institute of Fundamental Research,\\
		Homi Bhabha Rd, Bombay 400 005, INDIA.\\
e-mail address: orion@tifrvax.tifr.res.in, kirti@motive.math.tifr.res.in
			}

\end{document}